\documentclass[11pt,leqno]{article}
\usepackage{amsfonts}
\usepackage{amsmath}
\usepackage{amssymb}
\newcommand{\pref}[1]{(\ref{#1})}
\newtheorem{theo}{Theorem}[section]
\newtheorem{lem}[theo]{Lemma}

\newtheorem{defi}[theo]{Definition}

\title{\Large \bf {On weak holonomy}}
\author{{\sc Bogdan Alexandrov}
\thanks{Supported by SPP 1154 "Globale Differentialgeometrie" of DFG and The European Contract Human Potential Programme, Research Training Network HPRN-CT-2000-00101}
}
\date{}
\begin{document}
\maketitle
\vspace{5mm}
\begin{abstract}
We prove that $SU(n)$ ($n \ge 3$) and $Sp(n)U(1)$ ($n \ge 2$) are the only connected Lie groups acting transitively and effectively on some sphere which can be weak holonomy groups of a Riemannian manifold without having to contain its holonomy group. In both cases the manifold is K\"ahler.
\\[10mm]
{\bf Keywords:} holonomy group, weak holonomy group \\
{\bf MSC 2000: } 53C29
\\[10mm]
\end{abstract}

\section{Introduction}

Let $M$ be a Riemannian manifold and $p$ be a fixed point in $M$. Recall that the holonomy group $Hol(p)$ at $p$ is the group consisting of all parallel translations $\tau_\gamma$ along piecewise differentiable loops $\gamma$ based at $p$. The restricted holonomy group $Hol^0 (p)$ is the identity component  of $Hol(p)$ and consists of those $\tau_\gamma$ for which $\gamma$ is homotopic to 0.

In 1955 M. Berger \cite{Ber} found the finite list of the possible restricted holonomy groups of the Riemannian manifolds which are neither locally reducible nor locally symmetric. Later  Simons \cite{Sim} gave another proof of this result by showing that the hypothesis that the manifold is locally irreducible and not locally symmetric implies that $Hol(p)$ acts transitively on the unit sphere in the tangent space $T_p M$. The list of connected Lie groups acting effectively and transitively on some sphere is \cite{MS,Bor1,Bor2}
\begin{equation}\label{*}
SO(n), U(n), SU(n), Sp(n), Sp(n)U(1), Sp(n)Sp(1), G_2, Spin(7), Spin(9).
\end{equation}
For each of these groups its representation on $T_p M$ is isomorphic to its standard representation: $\mathbb{R}^n$ for $SO(n)$, $\mathbb{R}^{2n} \cong \mathbb{C}^n$ for $U(n)$ and $SU(n)$, $\mathbb{R}^{4n} \cong \mathbb{H}^n$ for $Sp(n)$, $Sp(n)U(1)$ and $Sp(n)Sp(1)$, $\mathbb{R}^7$ for $G_2$, $\mathbb{R}^8$ for $Spin(7)$, $\mathbb{R}^{16}$ for $Spin(9)$. Now it is an easy observation that $Sp(n)U(1)$ cannot be a (restricted) holonomy group and so one obtains Berger's list. It was proved later that $Hol^0 (p) = Spin(9)$ implies that the Riemannian manifold is locally symmetric \cite{Ale,Gr2} and that the other groups indeed appear as restricted holonomy groups of non locally symmetric spaces (see e.g. \cite{Bes,Br,J}).

It is well known that the Lie algebra of the holonomy group is determined by the curvature. Thus, if $Hol^0 (p) \not = SO(n)$, the curvature tensor will have additional properties. Motivated by the search for weaker conditions which will imply that the curvature satisfies some additional identities, A. Gray introduced in 1971 the notion of weak holonomy group \cite{Gr}. To define it one needs the auxiliary notion of special subspace for a group. We will give two definitions. 

\begin{defi}\label{def1}
{\rm Let $T$ be a real representation of the group $G$. A subspace $P \subseteq T$ is called {\it special} if the following conditions are satisfied:

(i) There exists a proper subspace $P' \subset P$ such that $g|_P$ is determined by $g|_{P'}$ for all $g \in G$.

(ii) If $P' \subset P \subseteq P''$ and $g|_{P''}$ is determined by $g|_{P'}$ for all $g \in G$, then $P'' = P$.}
\end{defi}

We call $P'$ a {\it generating subspace} for the special subspace $P$.

Let us give some examples of special subspaces.

{\bf 1)} For the standard representation of $SO(n)$ on $\mathbb{R}^n$ the only special subspace is $P = \mathbb{R}^n$ and each $(n-1)$-dimensional subspace $P'$ is a generating subspace.

{\bf 2)} For the standard representation of $U(n)$ on $\mathbb{R}^{2n} \cong \mathbb{C}^n$ the subspaces $P = span \{ x, Ix \}$ are special ($I$ is the $U(n)$-invariant complex structure) and $P' = span \{ x \}$ is a generating subspace for $P$.

{\bf 3)} Consider the standard representation of $G_2$ on $\mathbb{R}^7$. Let $\Phi : \mathbb{R}^7 \times \mathbb{R}^7 \longrightarrow \mathbb{R}^7$ be the $G_2$-invariant vector cross product \cite{Gr3}. Then each 2-dimensional subspace $P' = span \{ x, y \}$ generates a special subspace $P = span \{ x, y, \Phi (x,y) \}$.

\vspace{3mm}
In these three examples the dimension of the special subspaces is minimal.

\begin{defi}\label{def2}
{\rm Let $T$ be a real representation of the group $G$. A subspace $P \subseteq T$ is called {\it special} if the following conditions are satisfied:

(i) There exists a proper subspace $P' \subset P$ such that $g(P)$ is determined by $g(P')$ for all $g \in G$.

(ii) If $P' \subset P'' \subseteq P$ and $g(P'')$ is determined by $g(P')$ for all $g \in G$, then $P'' = P$.}
\end{defi}

Again we call $P'$ a {\it generating subspace} for the special subspace $P$. The above three  examples are also examples of special subspaces according to Definition~\ref{def2} which have minimal dimension.

The reason we give two different (and clearly non-equivalent in general) definitions of special subspaces is the following. In the original paper of A. Gray \cite{Gr} the definition of a special subspace is as Definition~\ref{def1} but in (ii) the restrictions $g|_{P''}$ and $g|_{P'}$ are replaced by the images $g(P'')$ and $g(P')$ respectively. In \cite{Gr2} the restrictions are replaced by images in both (i) and (ii). It is clear that there are typographical errors in these definitions since in this form they are useless: the only special subspace would be the whole space $T$. The most straightforward correction is Definition~\ref{def1} and this is also the definition used by L. Schwachh\"ofer in \cite{Schw}. In Definition~\ref{def2} we assume that A. Gray meant indeed images instead of restrictions but in (ii) $P''$ and $P$ were misplaced.

Let $G$ be a connected Lie subgroup of $SO(n)$ and $P_G$ be a reduction of the structure group $O(n)$ of $M$ to $G$, i.e., $P_G$ is a principal $G$-bundle over $M$ which is a subbundle of the bundle of orthonormal frames. Each element $u \in P_G$ defines an isometry between $\mathbb{R}^n$ and $T_p M$, where $p$ is the projection of $u$. Fixing such a $u$ allows us to consider $G$ as a subgroup of $SO(T_p M)$ and also $Hol(p)$ as a subgroup of $SO(n)$.

\begin{defi}\label{def3}
{\rm The group $G$ is called a {\it weak holonomy group} of $M$ if for each $p \in M$ and each differentiable loop $\gamma$ in $M$ with $\gamma (0) = \gamma (1) = p$ there exists $g \in G$ such that $\tau _\gamma |_P = g|_P$ whenever $P$ is a special subspace of $T_p M$ of minimal dimension with $\dot{\gamma} (0) \in P$.}
\end{defi}

Clearly, the definition does not depend on the choice of $u \in P_G$.

Notice that a Riemannian manifold may have more than one weak holonomy group. For example, any group $G$ such that $Hol(p) \subseteq G \subseteq SO(n)$ is a weak holonomy group. Of course, the interesting situation is when $G$ does not contain $Hol(p)$.

In \cite{Gr,Gr4} A. Gray studied the question which groups from the list \pref{*} can be weak holonomy groups of a Riemannian manifold without containing its holonomy group. He proved that the groups $U(n)$, $SU(n)$ and  $G_2$ have this property, the groups $SO(n$) (for trivial reasons), $Sp(n)$, $Sp(n)U(1)$, $Sp(n)Sp(1)$ and $Spin(7)$ do not, and the case $G=Spin(9)$ remained unresolved. According to his results, weak holonomy group $U(n)$ or $SU(n)$ is equivalent to the manifold being nearly K\"ahler, while weak holonomy group $G_2$ is equivalent to the existence of a nearly parallel vector cross product on $M$. His proofs relied on the following mysterious argument: If $S$ is certain $G$-invariant tensor on $M$ (for example, the complex structure $I$ if $G$ is $U(n)$ or $SU(n)$ or the vector cross product $\Phi$ if $G$ is $G_2$), then $G$ is a weak holonomy group if and only if $\nabla_X S(X,X_2,\dots,X_k) = 0$ for all $X,X_2,\dots,X_k$.

Unfortunately, it turns out that this argument is not correct. In section~\ref{sec2} we show that the simplest example of a nearly K\"ahler manifold, the 6-dimensional sphere $S^6$, does not have weak holonomy $U(3)$ and that the simplest example of a manifold with nearly parallel vector cross product, the 7-dimensional sphere $S^7$, does not have weak holonomy $G_2$. Notice that this does not depend on that which of \ref{def1} and \ref{def2} is the correct definition of a special subspace since the special subspaces for $U(n)$, $SU(n)$ and $G_2$ determined by Gray are the same as in the above examples.

Because of this we consider again the above mentioned question studied by Gray. We prove the following.

\begin{theo}\label{th1}
Whichever of Definition~\ref{def1} and Definition~\ref{def2} is used to define the notion of special subspace, the only connected Lie groups which act effectively and transitively on some  sphere that can be weak holonomy groups of a Riemannian manifold without containing its  holonomy group are $SU(n)$, $n \ge 3$, and $Sp(n)U(1)$, $n \ge 2$. In both cases the manifold is K\"ahler.
\end{theo}

Simple examples of manifolds with weak holonomy $SU(n)$ or $Sp(n)U(1)$ can be obtained by taking small open neighbourhoods in Riemannian manifolds with holonomy $U(n)$ or $U(2n)$ respectively. An interesting consequence of Theorem~\ref{th1} is that it is possible to have groups $G' \subset G'' \subset G'''$ so that $G'$ and $G'''$ are weak holonomy groups but $G''$ is not. For example, take $G' = Sp(n)U(1)$, $G'' = Sp(n)Sp(1)$, $G''' = SO(4n)$. Another curious fact is that, with the exception of $Sp(n)$, for all groups from the list \pref{*} the special subspaces according to Definition~\ref{def1} and according to Definition~\ref{def2} are the same.

\section{Examples}\label{sec2}

Let $<\cdot,\cdot>$ be the standard inner product in $\mathbb{R}^{n+1}$. Denote by  $e_1,\dots,e_{n+1}$ the standard basis of $\mathbb{R}^{n+1}$ and by $e^1,\dots,e^{n+1}$ its dual basis. We use the notation $e^{i_1 \dots i_k} := e^{i_1} \wedge \dots \wedge e^{i_k}$. Let $g_0$ be the standard metric on the unit sphere $S^n \subset \mathbb{R}^{n+1}$ and $\nabla$ (resp. $\nabla'$) be the Levi-Civita connection of $S^n$ (resp. $\mathbb{R}^{n+1}$). We have
$${\nabla'}_X Y = \nabla_X Y - <X,Y_x>x, \qquad X \in T_x S^n, \quad Y \in \Gamma(TS^n).$$
Let $x_t$ be a smooth curve on $S^n$ and let $\tau^0_s$ denote the parallel translation from $T_{x_0} S^n$ to $T_{x_s} S^n$ along $x_t$. If $X_0 \in T_{x_0} S^n$ and $X_t := \tau^0_t (X_0)$, then $X_t$ satisfies
$$\frac{d X_t}{d t} = - <\dot{x}_t,X_t>x_t.$$
Let $\omega$ be a parallel $(k+1)$-form on $\mathbb{R}^{n+1}$. Define $\Omega_x (X_1,\dots,X_k) = \omega (x,X_1,\dots,X_k)$ for $x \in S^n$. Then $\Omega$ is a $k$-form on $S^n$ which satisfies
$$\nabla_X \Omega (X,X_2,\dots,X_k) = 0 \qquad \mbox{for all } X,X_2,\dots,X_k.$$

\vspace{3mm}
\noindent {\bf Example 1} Consider the 2-form $\Omega$ on $S^6$, defined by $\Omega_x (X,Y) = \varphi (x,X,Y)$, where the 3-form $\varphi$ on $\mathbb{R}^7$ is
$$\varphi = e^{123} + e^{145} - e^{167} + e^{246} + e^{257} + e^{347} - e^{356}.$$
The 2-form $\Omega$ is non-degenerate and therefore it defines an almost complex structure $I$ compatible with $g_0$ on $S^6$ via the equation $g_0 (IX,Y) = \Omega (X,Y)$. The 3-form $\varphi$ is parallel on $\mathbb{R}^7$, so $\Omega$ satisfies $\nabla_X \Omega (X,Y) = 0$. This means, by definition, that $(S^6,g_0,I)$ is a nearly K\"ahler manifold.

Fix $r \in (0,1)$ and let
$$x_t = (r\cos t, r\sin t, 0,\sqrt{1-r^2}, 0,0,0) \in S^6.$$
Let $X_0 := \dot{x}_0 = (0,r,0,0,0,0,0)$. Then
$$X_t := \tau^0_t (X_0) = (- \frac{r^3 \sin((1+\sqrt{1-r^2})t)}{2(1+\sqrt{1-r^2})} - \frac{r^3 \sin((1-\sqrt{1-r^2})t)}{2(1-\sqrt{1-r^2})}, $$
$$\frac{r^3 \cos((1+\sqrt{1-r^2})t)}{2(1+\sqrt{1-r^2})} + \frac{r^3 \cos((1-\sqrt{1-r^2})t)}{2(1-\sqrt{1-r^2})},0, -r^2 \sin(\sqrt{1-r^2}t),0,0,0).$$
We have $Y_0 := IX_0 = (0,0,r^2,0,0,-r\sqrt{1-r^2},0)$ and $Y_t := \tau^0_t (Y_0) = Y_0$.

The subspace $P:= span \{ X_0,Y_0 \}$ is a special subspace of minimal dimension for $U(3)$ and $SU(3)$. We have $x_{2\pi} = x_0$ and $\dot{x}_0 = X_0 \in P$. If $U(3)$ (or $SU(3)$) were a weak holonomy group for $S^6$, then there would exist $g \in U(3)$ (or $SU(3)$) such that $\tau^0_{2\pi} |_P = g|_P$. This means
$$X_{2\pi} = \tau^0_{2\pi} (X_0) = g(X_0), \qquad Y_{2\pi} = \tau^0_{2\pi} (Y_0) = g(Y_0) = g(IX_0) = Ig(X_0).$$
Hence $IX_{2\pi} = Y_{2\pi} =Y_0 = IX_0$, i.e., $X_{2\pi} = X_0$. But $X_{2\pi}^2 = r \cos(\sqrt{1-r^2}.2\pi) \not = r = X_0^2$, a contradiction.

Thus, although $(S^6,g_0,I)$ is nearly K\"ahler, it does not have weak holonomy group $U(3)$ or $SU(3)$.

\vspace{3mm}
\noindent {\bf Example 2} Consider the 3-form $\varphi$ on $S^7$ defined by $\varphi_x (X,Y,Z) = \theta (x,X,Y,Z)$, where the 4-form $\theta$ on $\mathbb{R}^8$ is
$$\theta = e^{1234} + e^{1256} - e^{1278} + e^{1357} + e^{1368} + e^{1458} - e^{1467} $$
$$+ e^{5678} - e^{3456} + e^{2457} - e^{2358} + e^{2367} + e^{3478} + e^{2468}.$$
The 3-form $\varphi$ is non-degenerate and defines a reduction of the structure group of  $(S^7,g_0)$ to $G_2$. The corresponding vector cross product $\Phi$ is defined by $g_0 (\Phi(X,Y),Z) = \varphi (X,Y,Z)$. The 4-form $\theta$ is parallel on $\mathbb{R}^8$, so $\varphi$ satisfies $\nabla_X \varphi (X,Y,Z) = 0$. Thus $(S^7,g_0,\Phi)$ is a manifold with nearly parallel vector cross product \cite{Gr3}.

Fix $r \in (0,1)$ and let
$$x_t = (r\cos t, r\sin t, 0,\sqrt{1-r^2}, 0,0,0,0) \in S^7.$$
Let $X_0 := \dot{x}_0 = (0,r,0,0,0,0,0,0)$. Then
$$X_t := \tau^0_t (X_0) = (- \frac{r^3 \sin((1+\sqrt{1-r^2})t)}{2(1+\sqrt{1-r^2})} - \frac{r^3 \sin((1-\sqrt{1-r^2})t)}{2(1-\sqrt{1-r^2})}, $$
$$\frac{r^3 \cos((1+\sqrt{1-r^2})t)}{2(1+\sqrt{1-r^2})} + \frac{r^3 \cos((1-\sqrt{1-r^2})t)}{2(1-\sqrt{1-r^2})},0, -r^2 \sin(\sqrt{1-r^2}t),0,0,0,0).$$
Take $Y_0 := (0,0,r^2,0,0,-r\sqrt{1-r^2},0,0)$. Then $Y_t := \tau^0_t (Y_0) = Y_0$.

Let $Z_0 := \Phi(X_0,Y_0) = (-r^3 \sqrt{1-r^2},0,0,r^4,r^3 \sqrt{1-r^2},0,0,r^2 (1-r^2))$. Then
$$Z_t := \tau^0_t (Z_0) = (\frac{r^5 \cos((1+\sqrt{1-r^2})t)}{2(1+\sqrt{1-r^2})} - \frac{r^5 \cos((1-\sqrt{1-r^2})t)}{2(1-\sqrt{1-r^2})}, \frac{r^5 \sin((1+\sqrt{1-r^2})t)}{2(1+\sqrt{1-r^2})} $$
$$- \frac{r^5 \sin((1-\sqrt{1-r^2})t)}{2(1-\sqrt{1-r^2})},0, r^4 \cos(\sqrt{1-r^2}t),r^3 \sqrt{1-r^2},0,0,r^2 (1-r^2)).$$
The subspace $P:= span \{ X_0,Y_0,Z_0 \}$ is a special subspace of minimal dimension for $G_2$. We have $x_{2\pi} = x_0$ and $\dot{x}_0 = X_0 \in P$. If $G_2$ were a weak holonomy group for $S^7$, then there would exist $g \in G_2$ such that $\tau^0_{2\pi} |_P = g|_P$. This means
$$X_{2\pi} = \tau^0_{2\pi} (X_0) = g(X_0), \qquad Y_{2\pi} = \tau^0_{2\pi} (Y_0) = g(Y_0),$$ $$Z_{2\pi} = \tau^0_{2\pi} (Z_0) = g(Z_0) = g(\Phi(X_0,Y_0)) = \Phi(g(X_0),g(Y_0)),$$
i.e., $Z_{2\pi} = \Phi(X_{2\pi},Y_{2\pi})$. But
$$Z_{2\pi}^5 = r^3 \sqrt{1-r^2} \not = r^3 \sqrt{1-r^2} \cos(\sqrt{1-r^2}.2\pi) = \Phi(X_{2\pi},Y_{2\pi})^5,$$
a contradiction.

Thus, although $(S^7,g_0,\Phi)$ is a Riemannian manifold with nearly parallel vector cross product, it does not have weak holonomy group $G_2$.

\section{Determination of the minimal special subspaces}\label{sec3}

In this section we determine the special subspaces of minimal dimension for the connected Lie  groups acting transitively and effectively on some sphere according to Definition~\ref{def1} and Definition~\ref{def2} respectively.

\subsection{Minimal special subspaces according to Definition~\ref{def1}}\label{sec31}

Let $T$ be a real representation of the group $G$. For a subspace $P \subseteq T$ in this subsection we denote $G_P = \{ g \in G : g(x) = x \, \forall x \in P \}$. Clearly, if $P' \subseteq P$, then $G_{P'} \supseteq G_P$.

The following lemma is an obvious reformulation of Definition~\ref{def1}.

\begin{lem}\label{lem1}
Let $P' \subseteq T$ be a subspace and let $P$ be the maximal subspace of $T$ on which $G_{P'}$ acts trivially. If $P$ contains $P'$ strictly, then $P$ is a special subspace for $G$ generated by $P'$. Every special subspace arises in this way.
\end{lem}

If $P$ is a special subspace generated by $P'$, then $G_{P'} = G_P$. Thus two different subspaces $P_1'$ and $P_2'$ generate the same special subspace iff $G_{P_1'} = G_{P_2'}$.

As an extreme example let us consider the situation in which a special subspace $P$ is generated by $P' = \{ 0 \}$. Since $G_{P'} = G$, we obtain $G_P = G$, i.e., there exists a non-trivial subspace of $T$ on which $G$ acts trivially and $P$ is the maximal such subspace.

The groups we are interested in act transitively on the unit sphere in $T$. So for them $\dim P' \geq 1$ and therefore $\dim P \geq 2$. In particular, if we find a special subspace $P$ with $\dim P = 2$, then it is a minimal special subspace.

In the rest of this subsection we determine the minimal special subspaces for the groups from the list \pref{*}.

\vspace{3mm}
{\bf 1)} $G = SO(n)$

\nopagebreak
If $P' \subset \mathbb{R}^n$, $\dim P' = m$, then $G_{P'} \cong SO(n-m)$ and the orthogonal complement ${P'}^\bot$ of $P'$ is the standard real irreducible representation of $SO(n-m)$. Thus ${P'}^\bot$ contains a subspace on which $G_{P'}$ acts trivially only if $n-m=1$. Therefore the only special subspace for $SO(n)$ is $P = \mathbb{R}^n$, generated by an arbitrary $(n-1)$-dimensional subspace $P'$.

\vspace{3mm}
{\bf 2)} $G = U(n)$, $n \geq 2$

\nopagebreak
Let $I$ be the endomorphism of $\mathbb{R}^{2n}$ which corresponds to the multiplication by $i$ in $\mathbb{C}^n$ under the identification $\mathbb{C}^n \cong \mathbb{R}^{2n}$. Then, considered as a subgroup of $SO(2n)$,
$U(n) = \nolinebreak \{ g \in SO(2n): g(I) = I \}$.

Let $P' = span \{ x \} \subset \mathbb{R}^{2n}$. Then $G_{P'} \cong U(n-1)$ and
$$\mathbb{R}^{2n} = span \{ x, Ix \} \oplus span \{ x, Ix \}^\bot .$$
$G_{P'} \cong U(n-1)$ acts trivially on $span \{ x, Ix \}$ and $span \{ x, Ix \}^\bot$ is the standard real $2(n-1)$-dimensional irreducible representation of $U(n-1)$. Thus $P = span \{ x, Ix \}$ is a special subspace for $U(n)$ generated by $P'$ and has minimal dimension.

\vspace{3mm}
{\bf 3)} $G = SU(n)$, $n \geq 3$

\nopagebreak
Let $P' = span \{ x \} \subset \mathbb{R}^{2n}$. Then $G_{P'} \cong SU(n-1)$ and $\mathbb{R}^{2n}$ decomposes as in the previous case. $G_{P'} \cong SU(n-1)$ acts trivially on $span \{ x, Ix \}$, while $span \{ x, Ix \}^\bot$ is the standard real $2(n-1)$-dimensional representation of $SU(n-1)$, which is irreducible since $n-1 \geq 2$. Thus $P = span \{ x, Ix \}$ is a special subspace for $SU(n)$ generated by $P'$ and has minimal dimension.

\vspace{3mm}
{\bf 4)} $G = Sp(n)$

\nopagebreak
Let $I$, $J$, $K$ be the endomorphisms of $\mathbb{R}^{4n}$ which correspond respectively to the multiplication by $-i$, $-j$, $-k$ on the right in $\mathbb{H}^n$ under the identification $\mathbb{H}^n \cong \mathbb{R}^{4n}$. Then the group $Sp(n) = \{ A \in GL(n,\mathbb{H}): \bar{A}^t A = \pmb{1} \}$, considered as a subgroup of $SO(4n)$ under the above identification, is $Sp(n) = \{ g \in SO(4n): g(I) = I, g(J) = J, g(K) = K \}$.

Let $P' = span \{ x \} \subset \mathbb{R}^{4n}$. Then $G_{P'} \cong Sp(n-1)$ and
$$\mathbb{R}^{4n} = span \{ x, Ix, Jx, Kx \} \oplus span \{ x, Ix, Jx, Kx \}^\bot .$$
$G_{P'} \cong Sp(n-1)$ acts trivially on $span \{ x, Ix, Jx, Kx \}$, while $span \{ x, Ix, Jx, Kx \}^\bot$ is the standard real $4(n-1)$-dimensional irreducible representation of $Sp(n-1)$. Thus $P = span \{ x, Ix, Jx, Kx \}$ is a 4-dimensional special subspace generated by $P'$.

If $\dim P' >1$, then either $P' \subseteq span \{ x, Ix, Jx, Kx \}$ or there exists $0 \not = y \in P'$ such that $y \bot span \{ x, Ix, Jx, Kx \}$. In the first case the special subspace $P$ generated by $P'$ is the same as above. In the second one $G_{P'}$ acts trivially on $span \{ x, Ix, Jx, Kx \} \oplus span \{ y, Iy, Jy, Ky \}$, i.e., $\dim P \ge 8$.

Thus the minimal special subspaces for $Sp(n)$ are 4-dimensional and have the form $P = span \{ x, Ix, Jx, Kx \}$.

\vspace{3mm}
{\bf 5)} $G = Sp(n)U(1)$, $n \geq 2$

\nopagebreak
By definition $Sp(n)U(1) = (Sp(n) \times U(1))/_{\mathbb{Z}_2}$. We identify $U(1) \cong S^1 \subset \mathbb{C} \subset \mathbb{H}$. Then the action of $Sp(n)U(1)$ on $\mathbb{H}^n$ is given as follows: if $x \in \mathbb{H}^n$, $g=[A,a] \in Sp(n)U(1)$, where $A \in Sp(n)$, $a \in U(1)$, then $g(x) = Axa^{-1}$. Thus, under the identification $\mathbb{H}^n \cong \mathbb{R}^{4n}$, $Sp(n)U(1)$ becomes $Sp(n)U(1) = \{ g \in SO(4n): g(I) = I, g(span \{ J,K \}) = span \{ J,K \} \}$.

Let $P' = span \{ x \} \subset \mathbb{R}^{4n}$. Without loss of generality we can assume that $x=e_1 \in \mathbb{H}^n$. We have $G_{P'} \cong Sp(n-1)U(1)$, where the embedding $Sp(n-1)U(1) \hookrightarrow Sp(n)U(1)$ is given by $[B,a] \mapsto [A,a]$ with
$$A =
\begin{pmatrix}
a & 0 \\
0 & B
\end{pmatrix}
\in Sp(n) \subset GL(n,\mathbb{H}).$$
Now
$$\mathbb{R}^{4n} = span \{ x, Ix \} \oplus span \{ Jx, Kx \} \oplus span \{ x, Ix, Jx, Kx \}^\bot $$
and $G_{P'} \cong Sp(n-1)U(1)$ acts trivially on $span \{ x, Ix \}$, irreducibly on $span \{ Jx, Kx \}$ (because the $U(1)$-part acts irreducibly on this space) and $span \{ x, Ix, Jx, Kx \}^\bot$ is the standard real $4(n-1)$-dimensional irreducible representation of $Sp(n-1)U(1)$. Thus $P = span \{ x, Ix \}$ is a special subspace for $Sp(n)U(1)$ generated by $P'$ and has minimal dimension.

Notice that this result differs from the one in \cite{Gr} where it is claimed that the minimal special subspaces for $Sp(n)U(1)$ are 4-dimensional and have the form $span \{ x, Ix, Jx, Kx \}$.

\vspace{3mm}
{\bf 6)} $G = Sp(n)Sp(1)$, $n \geq 2$

\nopagebreak
By definition $Sp(n)Sp(1) = (Sp(n) \times Sp(1))/_{\mathbb{Z}_2}$. We identify $Sp(1) \cong S^3 \subset \mathbb{H}$. Then the action of $Sp(n)Sp(1)$ on $\mathbb{H}^n$ is given as  follows: if $x \in \mathbb{H}^n$, $g=[A,a] \in Sp(n)Sp(1)$, where $A \in Sp(n)$, $a \in Sp(1)$, then $g(x) = Axa^{-1}$. Thus, under the identification $\mathbb{H}^n \cong \mathbb{R}^{4n}$, $Sp(n)Sp(1)$ becomes $Sp(n)Sp(1) = \{ g \in SO(4n): g(span \{ I,J,K \}) = span \{ I,J,K \} \}$.

Let $P' = span \{ x \} \subset \mathbb{R}^{4n}$. Without loss of generality we can assume that $x=e_1 \in \mathbb{H}^n$. We have $G_{P'} \cong Sp(n-1)Sp(1)$, where the embedding $Sp(n-1)Sp(1) \hookrightarrow Sp(n)Sp(1)$ is given by $[B,a] \mapsto [A,a]$ with
$$A =
\begin{pmatrix}
a & 0 \\
0 & B
\end{pmatrix}
\in Sp(n) \subset GL(n,\mathbb{H}).$$
We have
$$\mathbb{R}^{4n} = P' \oplus span \{ Ix, Jx, Kx \} \oplus span \{ x, Ix, Jx, Kx \}^\bot .$$
$G_{P'} \cong Sp(n-1)Sp(1)$ acts irreducibly on $span \{ Ix, Jx, Kx \}$ (because this is isomorphic to the adjoint representation of the $Sp(1)$-part) and $span \{ x, Ix, Jx, Kx \}^\bot$ is the standard real $4(n-1)$-dimensional irreducible representation of $Sp(n-1)Sp(1)$. Thus the maximal subspace on which $G_{P'}$ acts trivially is $P'$ itself and therefore $P'$ does not generate any special subspace.

Let now $P' = span \{ x,y \}$ be 2-dimensional.

\vspace{3mm}
{\bf a)} $y \in span \{ Ix, Jx, Kx \}$

\nopagebreak
Without loss of generality we can assume that $y = Ix$. Then $G_{P'} \cong Sp(n-1)U(1)$ embedded in $Sp(n)U(1) \subset Sp(n)Sp(1)$ as in 5). As we saw in 5), $span \{ Jx, Kx \}$ and $span \{ x, Ix, Jx, Kx \}^\bot$ are irreducible representations of $Sp(n-1)U(1)$. Thus the maximal subspace on which $G_{P'}$ acts trivially is $P'$ itself and therefore $P'$ does not generate any special subspace.

\vspace{3mm}
{\bf b)} $y \in span \{ x, Ix, Jx, Kx \}^\bot$

\nopagebreak
Then $G_{P'} \cong Sp(n-2)Sp(1)$,
$$\mathbb{R}^{4n} = P' \oplus span \{ Ix, Jx, Kx \} \oplus span \{ Iy, Jy, Ky \} \oplus span \{ x, Ix, Jx, Kx, y, Iy, Jy, Ky \}^\bot $$
and $span \{ Ix, Jx, Kx \}$, $span \{ Iy, Jy, Ky \}$ and $span \{ x, Ix, Jx, Kx, y, Iy, Jy, Ky\}^\bot$ are irreducible \hfill representations \hfill of \hfill $Sp(n-2)Sp(1)$ \hfill (if \hfill $n=2$, \hfill then \hfill $G_{P'} \cong SO(3)$,  \\
$span \{ x, Ix, Jx, Kx, y, Iy, Jy, Ky\}^\bot = 0$ and $span \{ Ix, Jx, Kx \}$ and $span \{ Iy, Jy, Ky \}$ are isomorphic to the adjoint representation of $SO(3)$, i.e., again irreducible). So again $P'$ does not generate any special subspace.

\vspace{3mm}
{\bf c)} $y=y_1 + y_2$, $0 \not = y_1 \in span \{ Ix, Jx, Kx \}$, $0 \not = y_2 \in span \{ x, Ix, Jx, Kx \}^\bot$

\nopagebreak
Without loss of generality $x=e_1 \in \mathbb{H}^n$, $y_1 = \lambda Ix = \lambda Ie_1$, $\lambda >0$, $y_2 = e_2 \in \mathbb{H}^n$.

If $g \in Sp(n-1)Sp(1)$ is such that $g(y) = y$, then $g(y_1) = y_1$ and $g(y_2) = y_2$. Hence $G_{P'} \cong Sp(n-1)U(1) \cap Sp(n-2)Sp(1) = Sp(n-2)U(1)$, embedded in $Sp(n)Sp(1)$ through $[B,a] \mapsto [A,a]$ with
$$A =
\begin{pmatrix}
a & 0 & 0 \\
0 & a & 0 \\
0 & 0 & B
\end{pmatrix}
\in Sp(n) \subset GL(n,\mathbb{H}).$$
We have
$$\mathbb{R}^{4n} = span \{ x, Ix, y_2, Iy_2 \} \oplus span \{ Jx, Kx \} \oplus span \{ Jy_2, Ky_2 \} $$
$$\oplus span \{ x, Ix, Jx, Kx, y_2, Iy_2, Jy_2, Ky_2 \}^\bot .$$
$G_{P'}$ acts trivially on $span \{ x, Ix, y_2, Iy_2 \}$, while the other summands are irreducible representations of $G_{P'}$ (this is true also for $n=2$: in this case $G_{P'} = U(1)$ and the last summand is 0, but $\{ Jx, Kx \}$ and $span \{ Jy_2, Ky_2 \}$ are again irreducible representations of $U(1)$).

Thus the 4-dimensional subspaces $P = span \{ x, Lx, y_2, Ly_2 \}$ with $L \in span \{ I,J,K \}$, $y_2 \in span \{ x, Ix, Jx, Kx \}^\bot$, are special subspaces for $Sp(n)Sp(1)$ of minimal dimension.

If we take $P' = span \{ x, Ix, Jx \}$, then $G_{P'} \cong Sp(n-1)$ and
$$\mathbb{R}^{4n} = span \{ x, Ix, Jx, Kx \} \oplus span \{ x, Ix, Jx, Kx \}^\bot .$$
Thus we see as in 4) that $P = span \{ x, Ix, Jx, Kx \}$ is also a special subspace for $Sp(n)Sp(1)$ of minimal dimension.

It is straightforward to show that there are no other special subspaces for $Sp(n)Sp(1)$ of minimal dimension (but we shall not need that). Notice that in \cite{Gr} only the subspaces of the type $P = span \{ x, Ix, Jx, Kx \}$ were identified as minimal special subspaces for $Sp(n)Sp(1)$.

\vspace{3mm}
{\bf 7)} $G = G_2$

\nopagebreak
Let $\varphi$ be the 3-form on $\mathbb{R}^7$ defined in Example~1 in section~\ref{sec2} and the vector cross product $\Phi : \mathbb{R}^7 \times \mathbb{R}^7 \longrightarrow \mathbb{R}^7$ be defined by $<\Phi (x,y),z> = \varphi (x,y,z)$. Then, as a subgroup of $SO(7)$, $G_2 = \{ g \in SO(7) : g(\varphi) = \varphi \} = \{ g \in SO(7) : g(\Phi) = \Phi \}$.

Let $P' = span \{ x \} \subset \mathbb{R}^7$. Then $G_{P'} \cong SU(3)$. We have $\mathbb{R}^7 = P' \oplus {P'}^\bot $ and ${P'}^\bot $ is the standard real 6-dimensional irreducible representation of $SU(3)$. Therefore $P'$ does not generate any special subspace.

Let $P' = span \{ x,y \}$ be 2-dimensional. Then $G_{P'} \cong SU(2)$ (cf case 3) above) and $$\mathbb{R}^7 = span \{ x, y, \Phi (x,y) \} \oplus span \{ x, y, \Phi (x,y) \}^\bot .$$
$G_{P'} \cong SU(2)$ acts trivially on $span \{ x, y, \Phi (x,y) \}$, while $span \{ x, y, \Phi (x,y) \}^\bot $ is the standard real 4-dimensional irreducible representation of $SU(2)$. Thus $P = span \{ x, y, \Phi (x,y) \}$ is a special subspace for $G_2$ of minimal dimension (three) generated by $P'$.

\vspace{3mm}
{\bf 8)} $G = Spin(7)$

\nopagebreak
Let $\theta$ be the 4-form on $\mathbb{R}^8$ defined in Example~2 in section~\ref{sec2} and the triple vector cross product (cf \cite{Gr3}) $\Theta : \mathbb{R}^8 \times \mathbb{R}^8 \times \mathbb{R}^8 \longrightarrow \mathbb{R}^8$ be defined by $<\Theta (x,y,z),w> = \theta (x,y,z,w)$. Then, as a subgroup of $SO(8)$, $Spin(7) = \{ g \in SO(8) : g(\theta) = \theta \} = \{ g \in SO(8) : g(\Theta) = \nolinebreak \Theta \}$.

Let $P' = span \{ x \} \subset \mathbb{R}^8$. Then $G_{P'} \cong G_2$ and ${P'}^\bot $ is the standard real 7-dimensional irreducible representation of $G_2$. Hence $P'$ does not generate any special subspace.

Let $P' = span \{ x,y \}$ be 2-dimensional. Then $G_{P'} \cong SU(3)$ (cf case 7) above) and ${P'}^\bot $ is the standard real 6-dimensional irreducible representation of $SU(3)$. So again $P'$ does not generate any special subspace.

Let $P' = span \{ x,y,z \}$ be 3-dimensional. Then $G_{P'} \cong SU(2)$ and
$$\mathbb{R}^8 = span \{ x, y, z, \Theta (x,y,z) \} \oplus span \{ x, y, z, \Theta (x,y,z) \}^\bot .$$
$G_{P'} \cong SU(2)$ acts trivially on the first space, while the second is the standard real 4-dimensional irreducible representation of $SU(2)$. Thus $P = span \{ x, y, z, \Theta (x,y,z) \}$ is a special subspace for $Spin(7)$ of minimal dimension (four) generated by $P'$.

\vspace{3mm}
{\bf 9)} $G = Spin(9)$

\nopagebreak
This case resembles the case $G = Sp(n)Sp(1)$.

Let $I_1,\dots,I_9$ be the generators of the Clifford algebra $Cl(\mathbb{R}^9, -<\cdot,\cdot>)$, considered as endomorphisms of its 16-dimensional real representation $\mathbb{R}^{16}$. They satisfy
$$I_\alpha^2 = \pmb{1}, \qquad I_\alpha^* = I_\alpha, \qquad I_\alpha I_\beta = -I_\beta I_\alpha, \quad \alpha \not = \beta .$$
The group $Spin(9)$, considered as a subgroup of $SO(16)$, is \cite{Fr} $Spin(9) = \{ g \in SO(16) : g(span \{ I_1,\dots,I_9 \}) =  span \{ I_1,\dots,I_9 \} \}$.

Let $P' = span \{ x \} \subset \mathbb{R}^{16}$. Then $G_{P'} \cong Spin(7)$ and ${P'}^\bot = U \oplus V$, where $U$ is the real 7-dimensional irreducible representation of $Spin(7)$ (via the projection $Spin(7) \longrightarrow SO(7)$) and $V$ is the real 8-dimensional irreducible representation of $Spin(7)$ (i.e., the real spin representation). Hence $P'$ does not generate any special subspace.

Let $P' = span \{ x,y \}$ be 2-dimensional.

\vspace{3mm}
{\bf a)} $y \in U$

\nopagebreak
Then $G_{P'}$ is the subgroup of $Spin(7)$ which projects on $SO(6)$, i.e., $G_{P'} \cong Spin(6)$. We have $\mathbb{R}^{16} = P' \oplus W_1 \oplus V$, where $W_1$ is the orthogonal complement of $span \{ y \}$ in $U$. Hence $W_1$ is the real 6-dimensional irreducible representation of $Spin(6)$ (through the projection $Spin(6) \longrightarrow SO(6)$) and $V$ is the real 8-dimensional irreducible representation of $Spin(6)$ (i.e., the real spin representation). Thus in this case $P'$ does not generate any special subspace.

\vspace{3mm}
{\bf b)} $y \in V$

\nopagebreak
Then $G_{P'} \cong G_2$ (cf case 8) above). We have $\mathbb{R}^{16} = P' \oplus U \oplus W_2$, where $W_2$ is the orthogonal complement of $span \{ y \}$ in $V$. $U$ and $W_2$ are both isomorphic to the standard real 7-dimensional irreducible representation of $G_{P'} \cong G_2$. So again $P'$ does not generate any special subspace.

\vspace{3mm}
{\bf c)} $y=y_1 + y_2$, $0 \not = y_1 \in U$, $0 \not = y_2 \in V$

\nopagebreak
If $g \in Spin(7)$ is such that $g(y) = y$, then $g(y_1) = y_1$, $g(y_2) = y_2$. Hence $G_{P'} \cong Spin(6) \cap G_2 = SU(3)$ (because the subgroup of $G_2$ preserving $y_1 \in U$ is $SU(3)$).

Let $|x| = 1 = |y_2|$. Let $I \in span \{ I_1,\dots,I_9 \}$ be the endomorphism determined by $x$, i.e., the unique element of $span \{ I_1,\dots,I_9 \}$ such that $I^2 = \pmb{1}$, $Ix = x$ ($I$ is explicitly given by $I = \sum_{\alpha =1}^{9} <I_\alpha x,x>I_\alpha $). Then $span \{ x \} \oplus U$ and $V$ are the eigenspaces of $I$ for $+1$ and $-1$ respectively and $V = span \{ Jx : J \in span \{ I_1,\dots,I_9 \}, J \bot I \}$. Let $y_2 = Jx$. Then $J$ is $G_{P'}$-invariant. If $W_3$ is the orthogonal complement of $span \{ x,y_1 \}$ in $U$ and  $W_4$ the orthogonal complement of $span \{ y_2, Jy_1 \}$ in $V$, we have
$$\mathbb{R}^{16} = span \{ x, y_1, y_2, Jy_2 \} \oplus W_3 \oplus W_4 .$$
$G_{P'} \cong SU(3)$ acts trivially on $span \{ x, y_1, y_2, Jy_2 \}$, while $W_3$ and $W_4$ are isomorphic to the standard real 6-dimensional irreducible representation of $SU(3)$.

Thus $P = span \{ x, y_1, y_2, Jy_1 \}$ is a special subspace for $Spin(9)$ of minimal dimension generated by $P'$.

It is straightforward to see that every minimal special subspace for $Spin(9)$ has this form (but we shall not need that). Notice that the minimal special subspaces for $Spin(9)$ have dimension 4 in contrast to the result in \cite{Gr} where it is claimed that they are 8-dimensional.

\subsection{Minimal special subspaces according to Definition~\ref{def2}}\label{sec32}

Let $T$ be a real representation of the group $G$. For a subspace $P \subseteq T$ in this subsection we denote $G_P = \{ g \in G : g(P) = P \}$. The other notations are as in the previous subsection.

Definition~\ref{def2} is equivalent to the following.

\begin{lem}\label{lem2}
Let $P'$ be a proper subspace of $T$. Let $P$ be a subspace of $T$ which strictly contains $P'$ and such that $g(P) = P$ for each $g \in G_{P'}$. If there does not exist any subspace $P''$ such that $P' \subset P'' \subset P$ and $g(P'') = P''$ for each $g \in G_{P'}$, then  $P$ is a special subspace for $G$ generated by $P'$.
\end{lem}

Notice that, since $g(T) = T$ for each $g \in G_{P'}$, every proper subspace $P' \subset T$ generates a special subspace according to Definition~\ref{def2}, which is not the case with Definition~\ref{def1}.

The determination of the special subspaces simplifies for subgroups of $SO(n)$.

\begin{lem}\label{lem3}
Let $G$ preserve an inner product on $T$. Then $P \subseteq T$ is a special subspace generated by $P' \subset T$ iff $P = P' \oplus V$, where $\{ 0 \} \not = V \subseteq {P'}^{\bot}$ is an irreducible $G_{P'}$-representation.
\end{lem}

If $P' = \{ 0 \}$, then $G_{P'} = G$. Thus, if we suppose that $T$ is an irreducible representation of $G$, then the special subspace generated by $P' = \{ 0 \}$ is $T$. Therefore, when looking for the minimal special subspaces, we can always assume that the generating subspace $P'$ has $\dim P' \geq 1$ and hence the minimal special subspaces are at least 2-dimensional.

In the remainder of this subsection we determine the minimal special subspaces for the groups from the list \pref{*}.

\vspace{3mm}
{\bf 1)} $G = SO(n)$

\nopagebreak
If $P' \subset \mathbb{R}^n$, $\dim P' = m$, then $G_{P'} \cong SO(m) \times SO(n-m)$ and  ${P'}^\bot$ is irreducible with respect to $G_{P'}$. Thus the only special subspace for $SO(n)$ is the whole space $\mathbb{R}^n$.

\vspace{3mm}
{\bf 2)} $G = U(n)$, $n \geq 2$

\nopagebreak
Let $P' = span \{ x \} \subset \mathbb{R}^{2n}$. Then $G_{P'} \cong U(n-1)$ and the decomposition of ${P'}^{\bot}$ into irreducible $U(n-1)$-representations is
$${P'}^{\bot} = span \{ Ix \} \oplus span \{ x, Ix \}^\bot .$$
Thus the minimal special subspaces for $U(n)$ are $P = span \{ x, Ix \}$.

\vspace{3mm}
{\bf 3)} $G = SU(n)$, $n \geq 3$

\nopagebreak
Let $P' = span \{ x \} \subset \mathbb{R}^{2n}$. Then $G_{P'} \cong SU(n-1)$ and the decomposition of ${P'}^{\bot}$ into irreducible $SU(n-1)$-representations is the same as in the previous case. Hence the minimal special subspaces for $SU(n)$ are also $P = span \{ x, Ix \}$.

\vspace{3mm}
{\bf 4)} $G = Sp(n)$

\nopagebreak
Let $P' = span \{ x \} \subset \mathbb{R}^{4n}$. Then $G_{P'} \cong Sp(n-1)$ and we have the following  decomposition of ${P'}^{\bot}$ into irreducible $Sp(n-1)$-representations:
$${P'}^{\bot} = span \{ Ix \} \oplus span \{ Jx \} \oplus span \{ Kx \} \oplus span \{ x, Ix, Jx, Kx \}^\bot .$$
Therefore the minimal special subspaces for $Sp(n)$ are 2-dimensional and have the form $P = span \{ x, Lx \}$, where $L \in span \{ I,J,K \}$.

\vspace{3mm}
{\bf 5)} $G = Sp(n)U(1)$, $n \geq 2$

\nopagebreak
Let $P' = span \{ x \} \subset \mathbb{R}^{4n}$. Then $G_{P'} \cong Sp(n-1)U(1)$, embedded in  $Sp(n)U(1)$ as in 5) in the previous subsection. The decomposition of ${P'}^{\bot}$ into irreducible $G_{P'}$-representations is
$${P'}^{\bot} = span \{ Ix \} \oplus span \{ Jx, Kx \} \oplus span \{ x, Ix, Jx, Kx \}^\bot .$$
Hence the minimal special subspaces for $Sp(n)U(1)$ are 2-dimensional and have the form $P = span \{ x, Ix \}$.

\vspace{3mm}
{\bf 6)} $G = Sp(n)Sp(1)$, $n \geq 2$

\nopagebreak
Let $P' = span \{ x \} \subset \mathbb{R}^{4n}$. Then $G_{P'} \cong Sp(n-1)Sp(1)$, embedded in  $Sp(n)Sp(1)$ as in 6) in the previous subsection. The decomposition of ${P'}^{\bot}$ into irreducible $G_{P'}$-representations is
$${P'}^{\bot} = span \{ Ix, Jx, Kx \} \oplus span \{ x, Ix, Jx, Kx \}^\bot .$$
Thus \hfill the \hfill special \hfill subspace \hfill generated \hfill by \hfill $P'$, \hfill which \hfill has \hfill smaller \hfill dimension, \hfill is \\
$P = span \{ x, Ix, Jx, Kx \}$. We shall see that $P$ is indeed a special subspace of minimal dimension for $Sp(n)Sp(1)$.

Now let $P' = span \{ x,y \}$ be 2-dimensional.

\vspace{3mm}
{\bf a)} $y \in span \{ Ix, Jx, Kx \}$

\nopagebreak
Without loss of generality we can assume that $y = Ix$. Obviously $G_{P'} \supseteq Sp(n-1)U(1)$. The decomposition of ${P'}^{\bot}$ into irreducible $Sp(n-1)U(1)$-representations is
$${P'}^{\bot} =span \{ Jx, Kx \} \oplus span \{ x, Ix, Jx, Kx \}^\bot .$$ These spaces are also $G_{P'}$-invariant and therefore the special subspace of smaller dimension generated by $P'$ is again $P = span \{ x, Ix, Jx, Kx \}$.

\vspace{3mm}
{\bf b)} $y \in span \{ x, Ix, Jx, Kx \}^\bot$

\nopagebreak
Then $G_{P'} \supseteq Sp(n-2)Sp(1)$. The decomposition of ${P'}^{\bot}$ into irreducible $Sp(n-2)Sp(1)$-representations is
$${P'}^{\bot} = span \{ Ix, Jx, Kx \} \oplus span \{ Iy, Jy, Ky \} \oplus span \{ x, Ix, Jx, Kx, y, Iy, Jy, Ky \}^\bot $$
and the decomposition with respect to $G_{P'}$ cannot be finer. Therefore the special subspaces generated by $P'$ are at least 5-dimensional, i.e., they are not special subspaces  of minimal dimension for $Sp(n)Sp(1)$.

\vspace{3mm}
{\bf c)} $y=y_1 + y_2$, $0 \not = y_1 \in span \{ Ix, Jx, Kx \}$, $0 \not = y_2 \in span \{ x, Ix, Jx, Kx \}^\bot$

\nopagebreak
Without loss of generality $x=e_1 \in \mathbb{H}^n$, $y_1 = \lambda Ix = \lambda Ie_1$, $\lambda >0$, $y_2 = e_2 \in \mathbb{H}^n$. Then $G_{P'} \cong U(1) \times ((Sp(n-2) \times H)/_{\mathbb{Z}_2})$, where
$$H = \{ \alpha + \beta i : \alpha^2 + \beta^2 =1 \} \cup \{ \alpha j + \beta k : \alpha^2 + \beta^2 =1 \} \subset Sp(1)$$
(i.e.., $H \cong Pin(2)$). The group $U(1) \times ((Sp(n-2) \times H)/_{\mathbb{Z}_2})$ is embedded in $Sp(n)Sp(1)$ in the following way: $(e^{i\varphi},[B,a]) \mapsto [A,a] \in Sp(n)Sp(1)$, where
$$A =
\begin{pmatrix}
(\cos \varphi + \frac{i \lambda \sin \varphi}{\sqrt{1 + \lambda^2}}).a & -\frac{\varepsilon_a \sin \varphi}{\sqrt{1 + \lambda^2}}.a & 0\\
\frac{\sin \varphi}{\sqrt{1 + \lambda^2}}.a  & \varepsilon_a (\cos \varphi - \frac{i \lambda \sin \varphi}{\sqrt{1 + \lambda^2}}).a & 0\\
0 & 0 & B
\end{pmatrix}
\in Sp(n) \subset GL(n,\mathbb{H}),$$
with $\varepsilon_a = 1$ if $a = \alpha + \beta i$ and $\varepsilon_a = -1$ if $a = \alpha j + \beta k$ (i.e., $\varepsilon_a = \det (\pi(a))$, where $\pi : Pin(2) \longrightarrow O(2)$ is the projection). The decomposition of ${P'}^{\bot}$ into irreducible $G_{P'}$-representations is
\begin{eqnarray}
{P'}^{\bot} = span \{ -Ix + \lambda y_2, Iy_2 \} \oplus span \{ Jx, Kx, Jy_2, Ky_2 \} \nopagebreak[4] \nonumber \\
\oplus span \{ x, Ix, Jx, Kx, y_2, Iy_2, Jy_2, Ky_2 \}^\bot . \nonumber
\end{eqnarray}
Thus $P = span \{ x, Ix, y_2, Iy_2 \}$ is the special subspace of smallest dimension generated by $P'$.

So we see that the minimal special subspaces for $Sp(n)Sp(1)$ are 4-dimensional and that all subspaces $P = span \{ x, Lx, y, Ly \}$, where $L \in span \{ I,J,K \}$ and
$y \bot span \{ x, Lx \}$, are special for $Sp(n)Sp(1)$ of minimal dimension. It can be shown that these are all special subspaces for $Sp(n)Sp(1)$ of minimal dimension (but we shall not need that).

\vspace{3mm}
{\bf 7)} $G = G_2$

\nopagebreak
Let $P' = span \{ x \} \subset \mathbb{R}^7$. Then $G_{P'} \cong SU(3)$ and ${P'}^\bot $ is an irreducible representation of $SU(3)$. Thus the only special subspace generated by $P'$ is the whole space $\mathbb{R}^7$.

Let $P' = span \{ x,y \}$ be 2-dimensional. Then $G_{P'} \supseteq SU(2)$. The decomposition of ${P'}^{\bot}$ into irreducible $SU(2)$-representations is
$${P'}^{\bot} = span \{ \Phi (x,y) \} \oplus span \{ x, y, \Phi (x,y) \}^\bot .$$
These spaces are also $G_{P'}$-invariant and therefore this is also the decomposition of ${P'}^{\bot}$ into irreducible $G_{P'}$-representations. Thus the minimal special subspaces for $G_2$ are 3-dimensional and have the form $P = span \{ x, y, \Phi (x,y) \}$.

\vspace{3mm}
{\bf 8)} $G = Spin(7)$

\nopagebreak
Let $P' = span \{ x \} \subset \mathbb{R}^8$. Then $G_{P'} \cong G_2$ and ${P'}^\bot $ is an  irreducible representation of $G_2$. Thus the only special subspace generated by $P'$ is the whole space $\mathbb{R}^8$.

Let $P' = span \{ x,y \}$ be 2-dimensional. Then $G_{P'} \supseteq SU(3)$ and ${P'}^\bot $ is an irreducible representation of $SU(3)$ and therefore also of $G_{P'}$. So again the only special subspace generated by $P'$ is the whole space $\mathbb{R}^8$.

Let $P' = span \{ x,y,z \}$ be 3-dimensional. Then $G_{P'} \supseteq SU(2)$. The decomposition of ${P'}^{\bot}$ into irreducible $SU(2)$-representations is
$${P'}^{\bot} = span \{ \Theta (x,y,z) \} \oplus span \{ x, y, z, \Theta (x,y,z) \}^\bot .$$
These spaces are also preserved by $G_{P'}$ and therefore this is also the decomposition of ${P'}^{\bot}$ into irreducible $G_{P'}$-representations. Hence the minimal special subspaces for $Spin(7)$ are 4-dimensional and have the form $P = span \{ x, y, z, \Theta (x,y,z) \}$.

\vspace{3mm}
{\bf 9)} $G = Spin(9)$

\nopagebreak
Let $P' = span \{ x \} \subset \mathbb{R}^{16}$. Then $G_{P'} \cong Spin(7)$ and the decomposition of ${P'}^{\bot}$ into irreducible $Spin(7)$-representations is ${P'}^\bot = U \oplus V$. Thus the special subspaces generated by $P'$ are $P' \oplus U$ and $P' \oplus V$ and they have dimensions 8 and 9 respectively.

Let $P' = span \{ x,y \}$ be 2-dimensional.

\vspace{3mm}
{\bf a)} $y \in U$

\nopagebreak
Then $G_{P'} \supseteq Spin(6)$ and the decomposition of ${P'}^{\bot}$ into irreducible $Spin(6)$-re\-pre\-sen\-ta\-tions is ${P'}^{\bot} = W_1 \oplus V$. Thus the special  subspace of smallest dimension generated by $P'$ is at least 8-dimensional.

\vspace{3mm}
{\bf b)} $y \in V$

\nopagebreak
Then $G_{P'} \supseteq G_2$  and the decomposition of ${P'}^{\bot}$ into irreducible $G_2$-representations is ${P'}^{\bot} = U \oplus W_2$. Thus the special  subspace of smallest dimension generated by $P'$ is at least 9-dimensional.

\vspace{3mm}
{\bf c)} $y=y_1 + y_2$, $0 \not = y_1 \in U$, $0 \not = y_2 \in V$

\nopagebreak
Without loss of generality $|x| = 1 = |y_2|$. Let $|y_1| = \lambda$ and $I,J \in span \{ I_1,\dots,I_9 \}$ be as in 9c) in the previous subsection. We have $G_{P'} \supseteq SU(3)$ and the decomposition of ${P'}^{\bot}$ into irreducible $SU(3)$-representations is
$${P'}^{\bot} = span \{ -y_1 + \lambda^2 y_2 \} \oplus span \{ Jy_1 \} \oplus W_3 \oplus W_4 .$$
It is not hard to see that $span \{ -y_1 + \lambda^2 y_2, Jy_1 \}$ is $G_{P'}$-invariant and $G_{P'}$-irreducible. Thus the special subspace of smallest dimension generated by $P'$ is $P = span \{ x, y_1, y_2, Jy_1 \}$.

So the special subspaces of this type are special subspaces of minimal dimension for $Spin(9)$ (and every minimal special subspace for $Spin(9)$ has this form but we shall not use that).

\section{Proof of Theorem~\ref{th1}}

Let $G$ be a connected Lie subgroup of $SO(n)$ and $M$ be a Riemannian manifold whose structure group is reduced to $G$.

\begin{lem}\label{lem4}
Let the Lie group $G$ be compact. Suppose the subspace $P \subseteq T_p M$ is such that for each differentiable loop $\gamma$ in $M$ with $\gamma (0) = \gamma (1) = p$ and $\dot{\gamma} (0) \in P$ there exists $g \in G$ satisfying $\tau _\gamma |_P = g|_P$. Then for each $a \in Hol(p)$ there exists $g \in G$ such that $a|_P = g|_P$.
\end{lem}

\noindent {\it Proof:} Let $\gamma$ be a piecewise differentiable loop based at $p$ with $\dot{\gamma} (0) \in P$. There exists a continuous family of loops $\gamma_s$ based at $p$ such that $\gamma_0 = \gamma$, $\gamma_s$ is differentiable for $s>0$ and $\dot{\gamma}_s (0) = \dot{\gamma} (0)$. The parallel translation along a curve is a solution of an ODE and therefore $\tau _{\gamma_s}$ is continuous with respect to $s$. In particular, we obtain that there is a sequence of differentiable loops $\gamma_k$ such that $\dot{\gamma}_k (0) = \dot{\gamma} (0) \in P$ and $\tau _{\gamma_k} \longrightarrow \tau _\gamma$ as $k \longrightarrow \infty$. Thus there exist $g_k \in G$ such that $\tau _{\gamma_k} |_P = g_k |_P$. Because of the compactness of $G$ we can assume that the sequence $g_k$ is convergent. Let $g \in G$ be its limit. Then $\tau _\gamma |_P = g|_P$. Hence the property in the hypothesis of Lemma~\ref{lem4} is true also for piecewise differentiable loops.

Let $a \in Hol(p)$. Then $a = \tau _\delta$ for some piecewise differentiable loop $\delta$ based at $p$. If $\mu = \delta \circ \gamma ^{-1} \circ \gamma$, we have $\tau_\mu = \tau _\delta = a$. Since $\dot{\mu} (0) = \dot{\gamma} (0) \in P$, there exists $g \in G$ such that $\tau _\mu |_P = g|_P$, i.e., $a|_P = g|_P$.    \hfill $\Box$

\vspace{3mm}
Now let $M$ be a Riemannian manifold with weak holonomy $G$, where $G$ is one of the groups in the list \pref{*}. The results of the previous section show that if $G \not = Sp(n)$ the notion 'weak holonomy group $G$' is independent of which of the two definitions of special subspace is used.

\vspace{3mm}
{\bf 1)} $G = SO(n)$

\nopagebreak
There is nothing to prove in this case.

\vspace{3mm}
{\bf 2)} $G = U(n)$, $n \geq 2$

\nopagebreak
Let $a \in Hol(p)$. If $x \in T_p M$, then $P = span \{ x, Ix \}$ is a special subspace of minimal dimension. Definition~\ref{def3} and Lemma~\ref{lem4} imply that there exists $g \in U(n)$ such that $a(x) = g(x)$, $a(Ix) = g(Ix)$. But $g \circ I = I \circ g$. Hence $a(Ix) = Ia(x)$ for each $x \in T_p M$, i.e., $a \in U(n)$. Thus $Hol(p) \subseteq U(n)$.

\vspace{3mm}
{\bf 3)} $G = SU(n)$, $n \geq 3$

\nopagebreak
As in the previous case $Hol(p) \subseteq U(n)$. Manifolds with weak holonomy $SU(n)$ and holonomy $U(n)$ do exist. For example, every sufficiently small neighbourhood in a manifold with holonomy $U(n)$ has this property.

\vspace{3mm}
{\bf 4)} $G = Sp(n)$

\nopagebreak
Let $a \in Hol(p)$, $x \in T_p M$.

If we use Definition~\ref{def1} to define the special subspaces, then $P = span \{ x,Ix,Jx,Kx \}$ is a special subspace of minimal dimension. Hence, by Definition~\ref{def3} and Lemma~\ref{lem4} there exists $g \in Sp(n)$ such that
$$a(x) = g(x), \quad a(Ix) = g(Ix), \quad a(Jx) = g(Jx), \quad a(Kx) = g(Kx).$$
But $g$ commutes with $I$, $J$, $K$ and therefore
$$a(Ix) = Ia(x), \quad a(Jx) = Ja(x), \quad a(Kx) = Ka(x)$$
for each $x \in T_p M$. Thus also $a \in Sp(n)$, i.e., $Hol(p) \subseteq Sp(n)$.

If we use Definition~\ref{def2} to define the special subspaces, then $span \{ x, Ix \}$, $span \{ x, Jx \}$, $span \{ x, Kx \}$ are special subspaces of minimal dimension. As in case 2) this implies that $a$ commutes with $I$, $J$, $K$, i.e., again $Hol(p) \subseteq Sp(n)$.

\vspace{3mm}
{\bf 5)} $G = Sp(n)U(1)$, $n \geq 2$

\nopagebreak
As in case 2) $Hol(p) \subseteq U(2n)$. Manifolds with weak holonomy $Sp(n)U(1)$ and holonomy $U(2n)$ do exist. For example, one can take a sufficiently small neighbourhood in a manifold with holonomy $U(2n)$.

\vspace{3mm}
{\bf 6)} $G = Sp(n)Sp(1)$, $n \geq 2$

\nopagebreak
Let $a \in Hol(p)$, $x \in T_p M$, $L \in span \{ I,J,K \}$, $y \in span \{ x, Lx \}^\bot$.
Then $P = span \{ x, Lx, y, Ly \}$ is a special subspace of minimal dimension. By Definition~\ref{def3} and Lemma~\ref{lem4} there exists $g \in Sp(n)Sp(1)$ such that $a|_P = g|_P$. Since $g \in Sp(n)Sp(1)$, $g(L) = \alpha I + \beta J + \gamma K$. Hence
$$a(Lx) = g(Lx) = g(L)g(x) = g(L)a(x) = \alpha Ia(x) + \beta Ja(x) + \gamma Ka(x).$$
But $Ia(x)$, $Ja(x)$, $Ka(x)$ are linearly independent and so $\alpha$, $\beta$, $\gamma$ depend only on $a(Lx)$ and $a(x)$ and not on $g$.

The equality $a|_P = g|_P$ implies similarly
$$a(Ly) = \alpha Ia(y) + \beta Ja(y) + \gamma Ka(y)$$
with the same $\alpha$, $\beta$, $\gamma$, which depend only on $a(Lx)$ and $a(x)$). Thus this is true for each $y \in span \{ x, Lx \}^\bot$.

Finally, since $L^2 = -c.\pmb{1}$, we have
$$a(L^2 x) = -c.a(x) = -c.g(x) = g(L^2 x) = g(L)g(Lx) =  \alpha Ia(Lx) + \beta Ja(Lx) + \gamma Ka(Lx).$$
Thus $a(Lz) = \alpha Ia(z) + \beta Ja(z) + \gamma Ka(z)$ for each $z \in T_p M$.

Hence for each $L \in span \{ I,J,K \}$ there exist $\alpha$, $\beta$, $\gamma$ such that $a(L) = \alpha I + \beta J + \gamma K$. This means that $a \in Sp(n)Sp(1)$. Therefore $Hol(p) \subseteq Sp(n)Sp(1)$.

\vspace{3mm}
{\bf 7)} $G = G_2$

\nopagebreak
Let $a \in Hol(p)$ and $x, y \in T_p M$ be linearly independent. Then $P = span \{ x, y, \Phi (x,y) \}$ is a special subspace of minimal dimension. By Definition~\ref{def3} and Lemma~\ref{lem4} there exists $g \in G_2$ such that $a|_P = g|_P$. Since $g \in G_2$, $g(\Phi) = \Phi$. Hence
$$a(\Phi(x,y)) = g(\Phi(x,y)) = \Phi(g(x),g(y)) = \Phi(a(x),a(y)).$$
By continuity, $a(\Phi(x,y)) = \Phi(a(x),a(y))$ also when $x$ and $y$ are linearly dependent. Thus $a(\Phi) = \Phi$, i.e., $a \in G_2$. Therefore $Hol(p) \subseteq G_2$.

\vspace{3mm}
{\bf 8)} $G = Spin(7)$

\nopagebreak
In a similar way as in the previous case we see that $Hol(p) \subseteq Spin(7)$.

\vspace{3mm}
{\bf 9)} $G = Spin(9)$

\nopagebreak
This case is similar to case 6).

Let $a \in Hol(p)$ and $I, J \in span \{ I_1,\dots, I_9 \}$ be such that $I^2 = -\pmb{1}$, $I \bot J$. Let $x \in T_p M$ satisfy $Ix = x$ and $y_1 \in U$ (i.e., $Iy_1 = y_1$, $y_1 \bot x$). Then $P = span \{ x, y_1, Jx, Jy_1 \}$ is a special subspace of minimal dimension. By Definition~\ref{def3} and Lemma~\ref{lem4} there exists $g \in Spin(9)$ such that $a|_P = g|_P$. Since $g \in Spin(9)$, $g(I) = \sum_{\alpha =1}^9 \lambda_\alpha I_\alpha$.  Hence
$$a(Ix) = a(x) = g(x) = g(Ix) = g(I)g(x) = g(I)a(x) = \sum_{\alpha =1}^9 \lambda_\alpha I_\alpha a(x).$$
But $I_1 a(x),\dots, I_9 a(x)$ are linearly independent. So $\lambda_1,\dots, \lambda_9$ depend only on $a(x)$ and not on $g$.

The equality $a|_P = g|_P$ implies similarly
$a(Iy_1) = \sum_{\alpha =1}^9 \lambda_\alpha I_\alpha a(y_1)$
with the same $\lambda_1,\dots, \lambda_9$. Thus this is true for each $y_1 \in U$.

We have also
$$a(IJx) = -a(Jx) = -g(Jx) = g(IJx) = g(I)g(Jx) = g(I)a(Jx) = \sum_{\alpha =1}^9 \lambda_\alpha I_\alpha a(Jx)$$
and similarly $a(IJy_1) = \sum_{\alpha =1}^9 \lambda_\alpha I_\alpha a(Jy_1)$.

Thus $a(Iz) = \sum_{\alpha =1}^9 \lambda_\alpha I_\alpha a(z)$ for each $z \in T_p M$. Therefore $a(I) = \sum_{\alpha =1}^9 \lambda_\alpha I_\alpha $, i.e., $a(I) \in span \{ I_1,\dots, I_9 \}$ for each $I \in span \{ I_1,\dots, I_9 \}$. Hence $a \in Spin(9)$, i.e., $Hol(p) \subseteq Spin(9)$.

\vspace{3mm}
This completes the proof of Theorem~\ref{th1}.

\vspace{10mm}
\noindent
Bogdan Alexandrov \\
Universit\"at Greifswald \\
Institut f\"ur Mathemathik und Informatik \\
Friedrich-Ludwig-Jahn-Stra{\ss}e 15a \\
17487 Greifswald \\
Germany \\
{\tt e-mail: \quad boalexan@uni-greifswald.de}

\end{document}